\documentclass[11pt]{article}

\usepackage[a4paper,margin=1.15in]{geometry}
\usepackage{amsmath,amssymb,amsthm,mathtools}
\usepackage{enumitem}
\usepackage{microtype}
\usepackage{cite}
\usepackage[hidelinks]{hyperref}

\newtheorem{theorem}{Theorem}[section]
\newtheorem{lemma}[theorem]{Lemma}
\newtheorem{corollary}[theorem]{Corollary}
\newtheorem{proposition}[theorem]{Proposition}
\newtheorem{question}[theorem]{Question}
\theoremstyle{remark}
\newtheorem{remark}[theorem]{Remark}
\newtheorem{definition}[theorem]{Definition}

\theoremstyle{plain}
\newtheorem{maintheorem}{Theorem}

\newcommand{\R}{\mathbb R}
\newcommand{\Sph}{\mathbb S}
\newcommand{\Z}{\mathbb Z}
\newcommand{\Q}{\mathbb Q}

\newcommand{\Cech}{\check H}
\newcommand{\Hinf}{H_\infty}
\newcommand{\PD}{\operatorname{PD}}
\newcommand{\id}{\operatorname{id}}
\newcommand{\im}{\operatorname{im}}
\newcommand{\Hom}{\operatorname{Hom}}

\newcommand{\ord}{\operatorname{ord}}

\title{\bfseries Compact branch sets and injectivity of proper branched coverings}

\author{%
Deguang Zhong\\[0.4em]
\small Institute of Applied Mathematics, Shenzhen Polytechnic University,
Shenzhen 518055, China\\
\small \texttt{Email:huachengzhon@163.com}
}

\date{}

\begin{document}
\maketitle

\begin{abstract}
  In this paper, we show that a proper continuous open discrete map
  from the open unit ball $B^n$ onto a domain in $\R^n$ with $n\geq5$ is a homeomorphism whenever its branch set is compact. This answers Vuorinen's question affirmatively in every dimension $n\ge5.$
  In dimension four, we also show that, under the additional assumption of quasiregularity, the map is a quasiconformal homeomorphism.
  Together with the three-dimensional theorem of Kauranen, Luisto and Tengvall, these results answer the quasiregular question
  posed as Problem~7.66 in the collections of Barth--Brannan--Hayman and Hayman--Lingham in every dimension $n\ge3$.

\end{abstract}



\section{Introduction}\label{sec1}

Suppose that $B^n$ denote the open unit ball in $\R^n$.  In this paper, a continuous open
discrete map is called a branched covering \cite{R1993,KLT2022}. Its branch set $B_f$ consists of the points at which it is not a local
homeomorphism. We consider proper maps
\begin{equation}\label{eq1}
  f:B^n\longrightarrow Y=f(B^n)\subset\R^n .
\end{equation}
For properness we use the compact-preimage convention of \cite{Ho1975}. Its inverse images of compact subsets of the target $Y$ must be
compact in $B^n$ (Definition~\ref{def1}).

Following his work on exceptional sets and boundary behavior \cite{V1976}, Vuorinen \cite[Remarks~3.7]{V1979} asked whether,
for $n\ge3$, every proper branched covering $f$ as in \eqref{eq1} with compact branch set $B_f$ must be a homeomorphism onto $Y$. It was subsequently recorded in
\cite{V1981} and \cite[p.~193, problem~(4)]{V1988}. Srebro discussed the problem in his survey of the topology of quasiregular mappings
\cite[p.~108]{Sr1992}. Problem~7.66 in the collections of Barth, Brannan and Hayman \cite[p.~503]{BBH1984} and Hayman and Lingham
\cite[Problem~7.66]{HL2019} asks the corresponding question for a proper quasiregular map: does a compact branch set force injectivity in
dimension at least three? The term \emph{quasi-analytic} in the older formulation refers to the bounded distortion class used there. The
collections also record the positive special case $f(B^n)=B^n$ and the planar counterexample $z\mapsto z^2$. The original formulation already
concerns open discrete maps; it is not restricted to maps satisfying an analytic distortion inequality. The topological formulation of Vuorinen's question, with properness
understood relative to the image, reads as follows.

\begin{question}\label{que1}
  Let $n\ge3$ and let $f:B^n\to f(B^n)\subset\R^n$ be a proper continuous open discrete map with compact branch set. Must $f$ be a
  homeomorphism onto its image?
\end{question}
We first recall some work related to Question~\ref{que1}, including foundational results on branching, partial answers for maps from the Euclidean ball, and injectivity theorems in related settings.
 The local topology of open discrete mappings was developed by Church and Hemmingsen \cite{CH1960}, Chernavskii \cite{Ch1964}, and
V\"ais\"al\"a \cite{Va1966}. The Chernavskii--V\"ais\"al\"a theorem bounds the topological dimensions of the branch set and its image by $n-2$. V\"ais\"al\"a's degree theorem also shows that a proper open discrete map between oriented connected manifolds has finite, nonzero degree and constant multiplicity over the complement of the branch image. Rickman's monograph \cite{R1993} gives the local index, normal
neighborhood and path lifting theory in the form commonly used in geometric function theory. Fox's completed covering construction \cite{F1957} and the work of Berstein and Edmonds \cite{BE1978} provide related topological perspectives on branching.

Kauranen, Luisto and Tengvall proved the three-dimensional case of Vuorinen's question \cite[Theorem~1.1]{KLT2022}. They also proved the
empty-branch case in every dimension $n\ge2$ \cite[Theorem~1.2]{KLT2022}. Their main obstruction is torsion in the fundamental group at
infinity: if the target has torsion-free fundamental group at infinity, their Proposition~4.3 gives injectivity. In dimension three the
requisite torsion freeness follows from three-manifold topology \cite{P1957}. The same argument does not apply to higher-dimensional Euclidean
domains, whose fundamental groups can have torsion. Their paper left the general case in dimensions four and above unresolved. We use their
study of large saturated normal domains as the starting point for a more precise description of the covering at infinity.

The distinction between properness and weaker compactness conditions is essential in the history of this question. Apanasov \cite{A2019}
constructed bounded, locally homeomorphic quasiregular maps in three dimensions and discussed their relation to Vuorinen's injectivity
problem. Kauranen, Luisto and Tengvall \cite[Remark~4.6]{KLT2022} explained that these constructions are essentially proper, a condition which
does not imply properness relative to the image. In particular, restricting an essentially proper map to a slightly smaller ball does not in
general produce a proper map. Their empty-branch theorem excludes a proper noninjective example of this kind in every dimension.

There is a separate analytic history. Reshetnyak's theorem states that every nonconstant quasiregular map is open and discrete; see
\cite{Re1989,R1993}. Zorich's theorem \cite{Z1967} states that a branch-free entire quasiregular map $\R^n\to\R^n$, $n\ge3$, is a
quasiconformal homeomorphism. Martio, Rickman and V\"ais\"al\"a \cite{MRV1971} studied the topological and metric properties of quasiregular
maps, including isolated singularities. Kauranen, Luisto and Tengvall \cite{KLT2021} established compact-branch obstructions for entire
mappings of finite distortion under growth conditions on the distortion, and constructed examples showing the importance of those conditions.
In particular, the branch set of a nonconstant entire quasiregular map in dimension at least three is either empty or unbounded. These
entire-map results do not answer \eqref{eq1}: a source homeomorphism $B^n\cong\R^n$ preserves the topological hypotheses, but need not
preserve quasiregularity. For analytic background we refer also to \cite{IM2001,HK2014,V1988}.

Related injectivity results for proper maps were obtained by Bridson, Hinkkanen and Martin \cite{BHM2007}. Their argument uses the finite
index property of proper open surjections on fundamental groups, a consequence of the more general work of Smale \cite{Sma1957} and Walsh
\cite{Wa1976}. Their rigidity results include the statement that quasiregular self-maps of closed negatively curved Riemannian manifolds are
homeomorphisms when the dimension is not four; for closed hyperbolic manifolds the dimension restriction is unnecessary. These results concern
maps between manifolds with specified fundamental groups and curvature and do not by themselves answer the compact branch question for a
Euclidean ball. 

\subsection{Main results and proof strategy}
The following result answers Question \ref{que1} affirmatively in every dimension $n\ge5.$
\begin{maintheorem}\label{thm1}
  Let $n\ge5$ and let
  \[
    f:B^n\longrightarrow Y=f(B^n)\subset\R^n
  \]
  be a proper continuous open discrete map. If $B_f$ is compact, then $f$ is a homeomorphism.
\end{maintheorem}

In the four-dimensional setting, we use the analytic structure of a quasiregular map.
\begin{maintheorem}\label{thm2}
  Let $f:B^4\to Y=f(B^4)\subset\R^4$ be a proper nonconstant quasiregular map. If $B_f$ is compact, then $f$ is a quasiconformal homeomorphism
  onto $Y$.
\end{maintheorem}

 \emph{Our strategy to prove Theorem \ref{thm1} is topological.} No collar at infinity, finite homotopy type of an end neighborhood, or local flatness of the branch set is assumed.
After replacing the source by $\R^n$, we obtain connected cofinal end neighborhoods $Z_j$ whose full inverse images $C_j$ are regular covers
with a common finite deck group $G$. The group has order equal to the degree of $f$. The integral cohomology at infinity of the source is that
of $\Sph^{n-1}$. The resulting Cartan--Leray spectral sequence produces a full order class in $H^n(G;\Z)$.

Periodicity alone is not the decisive obstruction: finite periodic groups can have periods much smaller than $n$. The embedding
$Y\subset\Sph^n$ supplies a second piece of information. In degrees $1,\ldots,n-2$, cohomology at infinity splits additively into classes
extending over $Y$ and classes extending over its compact complement. A generator of a cyclic $p$-primary periodicity group must extend over
one of these sides. Its Bockstein preimage extends over the same side, as do their coefficient module cup products. A nonzero product in
degree $n-1$ would therefore extend over a side whose cohomology in that degree is zero. Swan's local period formulas \cite{Sw1960,N1974}
ensure that every nontrivial finite periodic group has a local period equal to two or four. This proves Theorem~\ref{thm1}.

\emph{In order to prove Theorem \ref{thm2}, we first note that the end deck transformations are quasiconformal.} Their extensions across the boundary sphere, together with local Smith theory, give a free
boundary action. The resulting quotient compactification $W$ is a compact locally quasiconformal four-manifold with boundary. Proper lifting
and the Euler characteristic of a finite covering show that $W$ is simply connected. Integral duality then shows that it is an integral
homology ball. If the degree exceeds one, its boundary is the Poincar\'e homology sphere. A boundary extension theorem of Tukia and
V\"ais\"al\"a \cite{TV1981} permits quasiconformal gluing to the $E_8$ plumbing. The resulting closed four-manifold is excluded by Donaldson
and Sullivan \cite{DS1989}.
\begin{remark}
The quasiregular version of Vuorinen's question is recorded as Problem~7.66 in Barth, Brannan and Hayman's
\emph{Research Problems in Complex Analysis} \cite[p.~503]{BBH1984} and in Hayman and Lingham's \emph{Research Problems in Function
Theory} \cite[p.~208]{HL2019}. Together with the three-dimensional result of Kauranen, Luisto and Tengvall
\cite[Theorem~1.1]{KLT2022}, Theorems~\ref{thm1} and~\ref{thm2} answer this problem affirmatively in every dimension $n\ge3$.
In dimension four, Theorem~\ref{thm2} addresses the quasiregular case of Vuorinen's question, using bounded distortion in both
the compactification and the gluing steps. It does not settle Vuorinen's general question for continuous open discrete maps in dimension four.
\end{remark}
\subsection{The connections between the technical tools for Theorems \ref{thm1} and \ref{thm2} and the existing literature} 
 Bounded length distortion mappings, introduced by Martio and V\"ais\"al\"a \cite{MV1988}, form another important subclass. Small-distortion
local invertibility results were obtained in \cite{KLT2021BLD}. Theorem~\ref{thm1} requires neither a length inequality nor Sobolev
regularity. This matters because branch sets need not admit smooth or piecewise linear models. Heinonen and Rickman \cite{HR1998} constructed
quasiregular maps with wild branching. Church and Timourian \cite{CT1978} studied small critical sets; suspension phenomena are discussed in
\cite{C1978}. Heinonen \cite{H2002} surveyed the limitations of the available branch set theory. Aaltonen \cite{A2016} studied monodromy
spaces of completed coverings; Aaltonen and Pankka \cite{AP2017} related local monodromy to branch set dimension. The piecewise linear
recognition theorem of Luisto and Prywes \cite{LP2021} has an additional hypothesis on the branch \emph{image}. We impose none of these
regularity conditions.

The group theoretic ingredients belong to the classical spherical space form theory of Milnor \cite{Mi1957}, Swan \cite{Sw1960R,Sw1960}, and
Thomas and Wall \cite{TW1971}. Our use of this theory concerns algebraic periodicity classes. It does not require an end neighborhood to have
finite CW type or to be an actual quotient of a sphere.

The deformation theory used in quasiconformal approximation and extension was developed by Edwards and Kirby \cite{EK1971}. Sullivan adapted
its construction to the Lipschitz and quasiconformal categories, and Tukia and V\"ais\"al\"a gave a detailed treatment, including manifolds
with boundary, in \cite{TV1981}. In dimensions at most three, Tukia \cite{Tu197879} had already obtained the Lipschitz Hauptvermutung by
piecewise linear approximation.

The analytic argument in dimension four uses a different part of quasiconformal theory. Sullivan \cite{Su1979} established the existence and
uniqueness of Lipschitz and quasiconformal manifold structures outside dimension four. Tukia and V\"ais\"al\"a \cite{TV1981} developed
approximation and extension results for manifolds with boundary. Their boundary extension theorem \cite[Theorem~4.10]{TV1981} applies in
dimension four as well: a locally quasiconformal boundary embedding extends to a locally quasiconformal embedding of a boundary neighborhood.
This theorem supplies the needed extension across the gluing boundary; a topological collar alone does not provide a quasiconformal product
collar. The extension theory from dimension $n$ to dimension $n+1$ was developed separately in \cite{TV1982}. Extension of quasiconformal maps
across spherical boundary components is treated by N\"akki \cite[Theorem~3.9]{Na1970}. For the round sphere arising in
Proposition~\ref{prop4}, we give a direct proof using Sobolev traces, reflection and Reshetnyak's theorem. Standard background on boundary
values and quasisymmetry is contained in \cite{Va1971,TV1980}.

Donaldson's theorem \cite{Do1983} restricts the definite intersection forms of closed simply connected smooth four-manifolds. Freedman's work
\cite{Fr1982} shows why this smooth conclusion cannot be imposed on an arbitrary topological filling. The relevant theorem here is the
quasiconformal obstruction of Donaldson and Sullivan \cite{DS1989}, whose analysis is carried out on four-manifolds with quasiconformal
coordinate changes. In particular, their argument in \cite[Section~7(ii)]{DS1989}, based on Fintushel and Stern \cite{FS1984}, excludes a
closed simply connected quasiconformal four-manifold with intersection form $-E_8$. This formulation avoids assuming that a quasiconformal
filling is smoothable. The plumbing description of the $E_8$ manifold and its Poincar\'e-sphere boundary is standard in four-manifold handle
theory; see \cite[Chapter~5]{GS1999}. The algebra of its intersection lattice is treated in \cite{MH1973}.

The boundary identification uses three-dimensional topology together with finite group cohomology. Moise's triangulation and uniqueness
theorems \cite{Mo1977} place topological three-manifolds in the usual piecewise linear and smooth categories, while Perelman's work
\cite{Pe2002,Pe2003RF,Pe2003} identifies closed simply connected three-manifolds and manifolds with finite fundamental group. The local fixed
set statements used for the boundary action belong to Smith theory, in the form developed in Borel's seminar \cite{Bo1960} and Bredon's
monograph \cite{Bd1972}. Compact manifolds with boundary are absolute neighborhood retracts; West's theorem \cite{We1977} gives them finite CW
homotopy type without requiring a triangulation of the four-dimensional filling.

There is also an extensive theory which represents four-manifolds by branched coverings. Montesinos \cite{Mo1978} represented oriented
four-dimensional $2$-handlebodies as three-sheeted covers of the ball branched over ribbon surfaces. Piergallini \cite{Pi1995} represented
closed oriented PL four-manifolds as four-sheeted covers of the sphere; Iori and Piergallini \cite{IP2002} removed the singularities of the
branch image after stabilization to five sheets. Piergallini and Zuddas \cite{PZ2019} extended this theory to compact manifolds with boundary
and open manifolds, including relative constructions. More recently, Bais, Piergallini and Zuddas \cite{BPZ2026} related the existence of
coverings between closed four-manifolds to embeddings of their intersection lattices, under a handle decomposition hypothesis on the target.
These representation results concern a different direction of construction and do not impose the compact branch end condition for a covering
whose source is an open Euclidean ball.

\section{Degree, duality and cohomology at infinity}\label{sec2}

In the rest of this paper, we assume that all manifolds are Hausdorff and second countable. They have no boundary unless a boundary is explicitly specified. Homology and ordinary
cohomology are singular; compact-support cohomology is denoted by $H_c^*$, and cohomology of arbitrary compact sets is \v{C}ech cohomology,
denoted by $\Cech^*$. We use constant coefficients unless an action is specified. The conventions for duality and coefficients are those of
\cite{Ha2002,Sp1966}. For the homology and homotopy terminology, see \cite[Chapters~0 and~2]{Ha2002}. An integrally \emph{acyclic} space has
zero reduced integral homology \cite{Ha2002}. The term \emph{rationally acyclic} has the analogous meaning with $\Q$ coefficients. A space is
\emph{contractible} if its identity map is homotopic to a constant map \cite[Chapter~0]{Ha2002}; this term will be used for the source in the
degree-transfer argument.

\subsection{Branched coverings and degree}
For branched coverings we use the terminology of \cite[Chapter~I]{R1993} and \cite[Section~2.2]{KLT2022}. For proper maps see \cite{Ho1975}.
The term essentially proper, used here for maps onto their target, comes from \cite{A2019}. See also \cite[Remark~4.6]{KLT2022}.

\begin{definition}\cite{R1993,KLT2022}\label{def1}
  A map $f:X\to Y$ is \emph{discrete} if every fiber is a discrete subset of $X$, and \emph{open} if it maps open sets to open sets. A
  \emph{branched covering} is a continuous open discrete map. Its \emph{branch set} is
  \[
    B_f=\{x\in X:f\text{ is not a local homeomorphism at }x\}.
  \]
  The set $f(B_f)$ is its branch image. The map is \emph{proper}  \cite{Ho1975} if $f^{-1}(K)$ is compact whenever $K\subset Y$ is compact. It
  is \emph{essentially proper}  \cite{A2019,KLT2022} if for every such $K$ there is a compact set $L\subset X$ with $f(L)=K.$
\end{definition}

For locally compact metrizable spaces, a continuous proper map is closed; see \cite{Ho1975}. A proper open map from a nonempty space to a
connected manifold is consequently surjective. Properness and discreteness make every fiber finite. A restriction $f:f^{-1}(V)\to V$ to an
open target set remains proper. If $V$ avoids the branch image, this restriction is a covering map \cite[Lemma~3.1]{KLT2022}.

For a proper branched covering between connected oriented $n$-manifolds, V\"ais\"al\"a's theorem \cite[Theorem~5.5]{Va1966} gives a finite
nonzero signed degree. Reversing one orientation if necessary, write this degree as $d>0$. The local indices $i(x,f)$ are positive integers
and satisfy
\begin{equation}\label{eq2}
  \sum_{x\in f^{-1}(y)}i(x,f)=d .
\end{equation}
At regular points the local index is one. Thus every fiber over $Y\setminus f(B_f)$ has exactly $d$ points. In particular, a proper branched
covering of degree one is a homeomorphism. The local index facts and normal neighborhood construction are treated in \cite[Chapter~I]{R1993}.

For normal domains and normal neighborhoods we follow \cite[Chapter~I, Section~4]{R1993} and \cite[Section~2.3]{KLT2022}.

\begin{definition}\cite[Section~2.3]{KLT2022}\label{def2}
  A domain $U$ compactly contained in the source is a \emph{normal domain} for $f$ if $\partial f(U)=f(\partial U)$. It is a
  \emph{normal neighborhood} of $x\in U$ if, in addition, $\overline U\cap f^{-1}(f(x))=\{x\}$.
\end{definition}

Here boundaries of source sets are taken in the source manifold. The restriction of a branched covering to a normal domain is proper onto its
image. Arbitrarily small normal neighborhoods exist \cite[Lemma~I.4.9]{R1993}.

For deck transformations and regular coverings, we follow \cite[Section~1.3]{Ha2002}. For a covering $p:E\to Z$, a \emph{deck transformation}
is a homeomorphism $g:E\to E$ satisfying $p\circ g=p$. These maps form the \emph{deck group} $\operatorname{Deck}(p)$ under composition
\cite{Ha2002}. On a connected covering a deck transformation is determined by its value at one point, so the deck group acts freely on every
fiber. A connected covering is \emph{regular} if this action is transitive on a fiber. Equivalently, $p_*\pi_1(E)$ is normal in $\pi_1(Z)$;
then the deck group is isomorphic to $\pi_1(Z)/p_*\pi_1(E)$. These statements use the usual path-connected, locally path-connected and
semilocally simply connected spaces, as holds for all covering spaces below \cite[Section~1.3]{Ha2002}.

We use the conventions on induced maps and duality in \cite[Sections~2.1, 3.1 and~3.3]{Ha2002}. For a continuous map $F:X\to N$, the induced
homomorphism $F_*:H_i(X;A)\to H_i(N;A)$ is obtained by composing singular simplices with $F$. Precomposing cochains with $F$ gives the
pullback $F^*:H^q(N;A)\to H^q(X;A)$. Properness is additionally needed for pullback on compact-support cohomology and pushforward on locally
finite homology. A locally finite chain has only finitely many simplices meeting each compact set. For an oriented $n$-manifold $N$, its
locally finite fundamental class is denoted by $[N]_{\mathrm{lf}}$. The \emph{cap product} \cite[Section~3.3]{Ha2002} evaluates a cochain on
part of a simplex and retains the complementary face, lowering the homological degree by the cohomological degree. In the form used here it is
\[
  H_c^q(N;\Z)\otimes H_n^{\mathrm{lf}}(N;\Z)
  \longrightarrow H_{n-q}(N;\Z),\qquad
  \alpha\otimes c\longmapsto\alpha\frown c.
\]
We write the cohomology factor first. Capping with $[N]_{\mathrm{lf}}$ is the Poincar\'e duality map $\PD_N$.

\begin{lemma}\label{lem1}
  Let $F:X\to N$ be a proper map of degree $d$ between connected oriented $n$-manifolds. There is an integral homomorphism
  \[
    F^!:H_i(N;\Z)\longrightarrow H_i(X;\Z)
    \quad\text{such that}\quad F_*F^!=d\,\id .
  \]
  If $X$ is contractible and $d\ne0$, then $d\,\widetilde H_i(N;\Z)=0$ in every degree.
\end{lemma}

\begin{proof}
  Noncompact Poincar\'e duality gives $\PD_N:H_c^{n-i}(N;\Z)\cong H_i(N;\Z)$ and the analogous isomorphism for $X$
  \cite[Theorem~3.35]{Ha2002}. Properness allows pullback in compact support cohomology. Set
  \[
    F^!=\PD_X\circ F^*\circ\PD_N^{-1}.
  \]
  Let $[X]_{\mathrm{lf}}$ and $[N]_{\mathrm{lf}}$ be the locally finite fundamental classes. For $\alpha\in H_c^{n-i}(N;\Z)$, cap-product
  naturality gives
  \[
    F_*\bigl(F^*\alpha\frown[X]_{\mathrm{lf}}\bigr)
    =\alpha\frown F_*[X]_{\mathrm{lf}}
    =d\,\alpha\frown[N]_{\mathrm{lf}}.
  \]
  This is the claimed identity. The relative fundamental classes and the cap-product projection formula are described in
  \cite[Section~3.3]{Ha2002}. The argument is integral and therefore applies to torsion classes as well.
\end{proof}

One can also obtain degree-annihilation identities from ramified covering transfers \cite{Sm1983,D1986}. Lemma~\ref{lem1} uses only proper
degree theory and manifold duality, so no identification with a particular definition of ramified covering is needed.

\begin{lemma}\label{lem2}
  Let $n\ge2$ and let $f:\R^n\to M$ be a proper branched covering of a connected orientable $n$-manifold. If $p:N\to M$ is the universal
  covering, then $f$ lifts to a proper surjective branched covering $F:\R^n\to N$. In particular, $\pi_1(M)$ is finite.
\end{lemma}

\begin{proof}
  The lifting criterion applies because $\R^n$ is simply connected \cite[Proposition~1.33]{Ha2002}. If $L\subset N$ is compact, then
  $F^{-1}(L)$ is closed in the compact set $f^{-1}(p(L))$. Thus $F$ is proper. Since $p$ is a local homeomorphism, $F$ is open and discrete.
  Its nonempty image is both open and closed in $N$, so $F$ is surjective.

  For $y\in M$, surjectivity implies $F(f^{-1}(y))=p^{-1}(y)$. The set on the left is finite. Hence the number of sheets of $p$, which equals
  $|\pi_1(M)|$, is finite.
\end{proof}

The finite index argument for proper open maps appears in \cite[Lemma~2.1]{BHM2007}, with antecedents in \cite{Sma1957,Wa1976}.
Lemma~\ref{lem2} records the properness of the lifted map as well, since the degree argument will be applied to that lift.

\subsection{Analytic subclasses}

We use the Sobolev definitions of quasiregular and quasiconformal mappings from \cite{R1993,IM2001}. For finite distortion and bounded length
distortion we use the conventions in \cite{HK2014} and the length formulation of \cite{MV1988,KLT2021BLD}, respectively.

A nonconstant continuous map $f:\Omega\to\R^n$ is \emph{quasiregular} if $f\in W_{\mathrm{loc}}^{1,n}(\Omega,\R^n)$ and
\[
  \|Df(x)\|^n\le KJ_f(x)\quad\text{for almost every }x
\]
for a finite constant $K\ge1$. Here $\|Df\|$ is the operator norm of the weak derivative and $J_f=\det Df$. A homeomorphic quasiregular map is
\emph{quasiconformal} \cite{R1993,IM2001}.

A map of \emph{finite distortion}  \cite{HK2014} allows a finite measurable function $K(x)$ in this inequality, with $f\in W_{\mathrm{loc}}^{1,1}$ and
$J_f\in L_{\mathrm{loc}}^1.$ A branched covering is \emph{$L$-BLD} \cite{MV1988,KLT2021BLD} if every path $\gamma$ satisfies,
with length allowed to be infinite,
\[
  L^{-1}\ell(\gamma)\le\ell(f\circ\gamma)\le L\ell(\gamma).
\]

\subsection{Ends and the compact-support sequence}

We use the neighborhood direct limit for cohomology at infinity, as in \cite[Section~7]{BC2025}. For the related compact-support construction
see \cite[Section~3.3, Exercises~17--18]{Ha2002} and \cite[Chapter~6]{Sp1966}.

\begin{definition}\cite[Section~7]{BC2025}\label{def3}
  For a locally compact space $X$ and an abelian group $A$, its \emph{cohomology at infinity} is
  \[
    \Hinf^q(X;A)=\varinjlim_{K\Subset X}H^q(X\setminus K;A).
  \]
  The limit runs over compact subsets, with transition maps induced by restrictions. A family of compact sets is \emph{cofinal} if every
  compact subset is contained in one member. Complements of a cofinal family are called cofinal neighborhoods of infinity.
\end{definition}

For inverse-limit terminology, see \cite[Section~3.F]{Ha2002}. The notation $\varprojlim^{\,1}$ denotes the first right derived functor of
inverse limit \cite{Ha2002}. Such a term can occur in inverse-limit calculations. No $\varprojlim^{\,1}$ term enters the filtered direct-limit
sequences below. We use the component description of ends in \cite[Section~3]{BC2025}. For a connected noncompact manifold, an \emph{end}
records a compatible choice of components outside compact subsets \cite{BC2025}. We shall construct a connected cofinal system directly.

\begin{lemma}\label{lem3}
  For a manifold $X$ there is a natural exact sequence
  \begin{equation}\label{eq3}
    \cdots\longrightarrow H_c^q(X;A)\longrightarrow H^q(X;A)
    \longrightarrow\Hinf^q(X;A)
    \longrightarrow H_c^{q+1}(X;A)\longrightarrow\cdots .
  \end{equation}
  For every space in this setting, $\Hinf^1(X;\Z)$ is torsion free.
\end{lemma}

\begin{proof}
  Take the filtered direct limit of the cohomology sequences of $(X,X\setminus K)$. Filtered direct limits of abelian groups are exact, and
  $H_c^q(X;A)=\varinjlim_K H^q(X,X\setminus K;A)$. This proves \eqref{eq3}. The universal coefficient theorem \cite{Ha2002} gives
  \[
    H^1(X\setminus K;\Z)
    \cong\Hom(H_1(X\setminus K;\Z),\Z),
  \]
  since $H_0$ is free abelian. Each group is torsion free. Multiplication by a nonzero integer is injective at every stage, and remains
  injective after an exact filtered direct limit. Thus the limit is torsion free. See also \cite[Section~3.3, Exercises~17--18]{Ha2002}.
\end{proof}

Euclidean exteriors retract onto spheres, so
\begin{equation}\label{eq4}
  \Hinf^q(\R^n;A)=
  \begin{cases}
    A,&q=0,n-1,\\
    0,&\text{otherwise}.
  \end{cases}
\end{equation}
For $n\ge3$, inclusions of sufficiently deep Euclidean exteriors kill their fundamental groups.

\subsection{The compact complement}

From now through Section~\ref{sec6}, let $f:\R^n\to Y\subset\R^n$ be a proper branched covering of positive degree $d$, with $n\ge3$. The
source change $B^n\cong\R^n$ is purely topological. View $Y$ as a subset of $\Sph^n=\R^n\cup\{\infty\}$ and put
\begin{equation}\label{eq5}
  D=\Sph^n\setminus Y .
\end{equation}
The set $D$ is nonempty and compact. Alexander duality takes the form
\begin{equation}\label{eq6}
  \widetilde H_i(Y;A)
  \cong\widetilde{\Cech}^{\,n-i-1}(D;A).
\end{equation}
The \v{C}ech version is essential because $D$ need not be locally contractible; see \cite[Section~3.3]{Ha2002} and \cite[Chapter~6]{Sp1966}.

\begin{lemma}\label{lem4}
  In the preceding setting,
  \[
    dH_i(Y;\Z)=0\quad(i>0),\qquad
    H_{n-1}(Y;\Z)=H_{n-2}(Y;\Z)=0.
  \]
  The compact set $D$ is connected. Moreover,
  \[
    H^q(Y;\Q)=\Cech^q(D;\Q)=0\quad(q>0).
  \]
  For every abelian coefficient group $A$,
  \begin{equation}\label{eq7}
    H^{n-1}(Y;A)=\Cech^{n-1}(D;A)=0,\qquad
    \Hinf^{n-1}(Y;A)\cong A,\quad \Hinf^n(Y;A)=0 .
  \end{equation}
\end{lemma}

\begin{proof}
  Lemma~\ref{lem1} annihilates all positive-dimensional homology of $Y$. The groups $\widetilde{\Cech}^0(D;\Z)$ and $\Cech^1(D;\Z)$ are
  torsion free. For degree zero this follows from locally constant integer-valued functions. For degree one it follows by taking the direct
  limit over polyhedral neighborhood approximations: each integer $H^1$ is torsion free, and filtered direct limits preserve this property.
  Alexander duality now proves the two integral vanishing statements. The degree-zero statement also makes $D$ connected.

  Tensoring homology with $\Q$ shows that $Y$ is rationally acyclic. The universal coefficient theorem with the injective coefficient group
  $\Q$ gives the asserted rational cohomology of $Y$. Alexander duality gives the rational cohomology of $D$ below degree $n$; in degrees at
  least $n$ it vanishes by continuity over proper open manifold neighborhoods of $D$ in $\Sph^n$. In top degree, connectedness of $Y$ implies
  $\Cech^{n-1}(D;A)=\widetilde H_0(Y;A)=0$. The integral vanishing in dimensions $n-1$ and $n-2$, followed by the universal coefficient
  theorem, gives $H^{n-1}(Y;A)=0$.

  A connected noncompact $n$-manifold has $H_i(Y;\Z)=0$ for $i\ge n$ \cite[Proposition~3.29]{Ha2002}; consequently $H^n(Y;A)=0$ here.
  Poincar\'e duality gives $H_c^n(Y;A)\cong H_0(Y;A)\cong A$ and $H_c^{n+1}(Y;A)=0$. Sequence~\eqref{eq3} proves the remaining statements in
  \eqref{eq7}.
\end{proof}

\section{Regularity at infinity}\label{sec3}
We now assume that $B_f$ is compact. We adapt the construction of normal domains in \cite[Lemma~4.1]{KLT2022}.  For any subset $A$ of the domain of $f$,  the set $f^{-1}(f(A))$ is used to
pass to complete fibers. We call this set the \emph{saturation}  \cite{KLT2022} of $A$ under $f.$ It is the smallest union of complete fibers
containing $A$. Properness makes the saturation $f^{-1}(f(B_f))$ compact as well. The following construction is a variant of
\cite[Lemma~4.1]{KLT2022}; it also makes explicit that the chosen end neighborhoods contain complete fibers.

\begin{lemma}\label{lem5}
  There are nested connected cofinal neighborhoods $Z_s$ of infinity in $Y$ such that $C_s=f^{-1}(Z_s)$ is connected, $f:C_s\to Z_s$ is a
  $d$-sheet covering, and, for suitable $R(s)>s$,
  \begin{equation}\label{eq8}
    \R^n\setminus\overline B_{R(s)}
    \subset C_s\subset\R^n\setminus\overline B_s .
  \end{equation}
  Both complements $Y\setminus Z_s$ and $\R^n\setminus C_s$ are compact.
\end{lemma}

\begin{proof}
  For sufficiently large $s$, let $K_s=f(\overline B_s)$ and $L_s=f^{-1}(K_s)$. These sets are compact, and $K_s$ contains the branch image.
  The map $\R^n\setminus L_s\to Y\setminus K_s$ is therefore a $d$-sheet covering. Let $C_s$ be the unique unbounded component of
  $\R^n\setminus L_s$. Uniqueness follows because $L_s$ is compact and a sufficiently large Euclidean exterior is connected. The image
  $Z_s=f(C_s)$ is a component of $Y\setminus K_s$: a component of the source of a covering maps onto a target component.

  The set $Z_s$ is not relatively compact in $Y$. Otherwise $C_s$ would lie in the compact inverse image of $\overline Z_s$. Every component
  of $f^{-1}(Z_s)$ maps onto $Z_s$. A component other than $C_s$ is bounded and has compact closure in $\R^n$; its image would make $Z_s$
  relatively compact in $Y$, a contradiction. Hence $f^{-1}(Z_s)=C_s$.

  Choose $R(s)$ with $L_s\subset B_{R(s)}$. Since $\overline B_s\subset L_s$, \eqref{eq8} follows. All source points outside $C_s$ lie in
  $\overline B_{R(s)}$. Surjectivity gives $Y\setminus Z_s\subset f(\overline B_{R(s)})$. The set on the left is closed in $Y$, so it is
  compact. The source complement is closed and bounded, hence compact. As $s$ increases, the $K_s$ and $L_s$ increase. The unbounded component
  for a larger parameter lies in the unbounded component for a smaller one, so the $C_s$ and $Z_s$ are nested. If $K\subset Y$ is compact,
  properness gives $f^{-1}(K)\subset\overline B_s$ for some $s$; surjectivity then gives $K\subset K_s\subset Y\setminus Z_s$. Likewise every
  compact source set lies in some $\overline B_s\subset\R^n\setminus C_s$. This proves cofinality on both sides.
\end{proof}

\begin{proposition}\label{prop1}
  After passing to a cofinal sequence of the neighborhoods in Lemma~\ref{lem5}, the coverings
  \[
    f:C_j\longrightarrow Z_j
  \]
  are regular with mutually compatible deck groups. These groups identify with one finite group $G$ of order $d$. Its action on
  $\Hinf^{n-1}(\R^n;\Z)\cong\Z$ is trivial.
\end{proposition}

\begin{proof}
  Fix $s$, and choose $t>R(s)$. By \eqref{eq8}, the inclusion $C_t\hookrightarrow C_s$ factors through $\R^n\setminus\overline B_{R(s)}$,
  which is simply connected. Thus
  \begin{equation}\label{eq9}
    \pi_1(C_t)\longrightarrow\pi_1(C_s)
    \quad\text{is trivial}.
  \end{equation}
  Choose a base point $y\in Z_t$ and a point $x$ of its full $d$-point fiber. Write $P_t=\pi_1(Z_t,y)$ and $j:P_t\to\pi_1(Z_s,y)$ for the
  inclusion homomorphism. Let $H_t\subset P_t$ be the subgroup for the covering $C_t\to Z_t$, based at $x$.

  If $\gamma\in H_t$, its lift based at $x$ is closed in $C_t$. By \eqref{eq9}, this lift is null-homotopic in $C_s$; hence $j(\gamma)=1$.
  Conversely, if $j(\gamma)=1$, the lift in the covering $C_s\to Z_s$ is closed. Since $C_t$ is the full preimage of $Z_t$, it is also the
  lift in $C_t$. Therefore
  \[
    H_t=\ker j .
  \]
  Because the lift criterion identifies $H_t$ with the stabilizer of $x$ in the monodromy action, the equality $H_t=\ker j$ makes this
  subgroup normal. The connectedness of $C_t$ makes the action on the full fiber transitive. This proves regularity, and $P_t/H_t\cong\im j$
  has order $[P_t:H_t]=d$. Equivalently, $\im j$ acts freely and transitively on the fiber.

  Take an increasing sequence $s_j$ with $s_{j+1}>R(s_j)$, and discard its first member. Every resulting cover is regular. A deck
  transformation of $C_j\to Z_j$ preserves $C_{j+1}=f^{-1}(Z_{j+1})$, so restriction gives a homomorphism of deck groups. It is injective
  because a deck transformation that is the identity on the nonempty open set $C_{j+1}$ is the identity on the connected cover $C_j$. Both
  groups have order $d$, so restriction is an isomorphism. Fixing these identifications gives $G$.

  Finally, $f\circ g=f$ for a deck transformation $g$. Locally both restrictions of $f$ are orientation-preserving homeomorphisms, so $g$
  preserves orientation. It also preserves every sufficiently deep saturated neighborhood $C_k$. Consequently it extends, fixing $\infty$, to
  a homeomorphism of $C_j\cup\{\infty\}$, an open neighborhood of $\infty$ in $\Sph^n$. This extension preserves orientation, and therefore
  acts by $+1$ on the local orientation class at $\infty$. The associated end class in $\Hinf^{n-1}(\R^n;\Z)$ is preserved as well.
\end{proof}

\begin{remark}\label{rem1}
  Proposition~\ref{prop1} is a statement about the end. It does not extend the deck transformations over the compact part of the source. Also,
  an individual $C_j$ need not be simply connected. Equation~\eqref{eq9} concerns the inclusion homomorphism, not the fundamental group of
  each neighborhood. The sandwich \eqref{eq8} gives
  \[
    \varinjlim_jH^q(C_j;A)\cong H^q(\Sph^{n-1};A),
  \]
  which is exactly the cohomological information needed below.
\end{remark}

\section{The end spectral sequence and periodicity}\label{sec4}

\subsection{Cartan--Leray at infinity}

For group modules, group cohomology and Tate cohomology we use the conventions of \cite[Chapters~I, III and~VI]{Br1982}. A $G$-module is an
abelian group with a left action of $G$ by automorphisms, equivalently a left module over the group ring $\Z G$. For a finite group $G$ and a
$G$-module $A$, write $H^*(G;A)$ for group cohomology and $\widehat H^*(G;A)$ for Tate cohomology. Here $H^q(G;A)$ is computed by applying
$\Hom_{\Z G}(-,A)$ to a projective resolution of the trivial module $\Z$. Tate cohomology extends this construction to all integral degrees by
a complete resolution. The Cartan--Leray spectral sequence relates the cohomology of a free group quotient to group cohomology with
coefficients in the cohomology of the covering space. More explicitly, put $A^G=\{a:ga=a\text{ for all }g\in G\}$, $N_Ga=\sum_{g\in G}ga$, and
let $I_GA$ be the subgroup generated by the elements $ga-a$. The low and negative Tate groups are
\[
  \widehat H^0(G;A)=A^G/N_GA,\qquad
  \widehat H^{-1}(G;A)=\ker(N_G)/I_GA,
\]
and $\widehat H^q(G;A)=H_{-q-1}(G;A)$ for $q\le-2$. For $q>0$, Tate and ordinary group cohomology agree \cite[Chapter~VI]{Br1982}.

The cup products below are the products with paired coefficients of \cite[Chapters~III, V and~VI]{Br1982} and \cite[Section~3.2]{Ha2002}. For
a commutative ring $R$ and an $RG$-module $M$, the $G$-equivariant pairing $M\otimes_R R\to M$, $m\otimes r\mapsto rm$, induces the
\emph{coefficient module cup product} \cite{Br1982,Ha2002}
\[
  H^i(G;M)\otimes_R H^j(G;R)\longrightarrow H^{i+j}(G;M).
\]
The same coefficient pairing defines cup products for spaces, neighborhood limits and Tate cohomology. Only the second factor requires a
coefficient ring; no multiplication on $M$ is needed. When a commutative ring is used as a coefficient group, its $G$-action will be trivial.
Cup products with a coefficient module use the module structure, rather than an unspecified multiplication on that module.

We use the cohomological indexing of \cite[Chapter~1, Theorem~1.14]{HaSS}. A cohomological spectral sequence consists of bigraded groups
$E_r^{a,b}$, called the terms on the $r$th page, and differentials \cite{HaSS}
\[
  d_r:E_r^{a,b}\longrightarrow E_r^{a+r,b-r+1}.
\]
The next page is the cohomology of this differential: $E_{r+1}^{a,b}=\ker(d_r|_{E_r^{a,b}})/ \im(d_r|_{E_r^{a-r,b+r-1}})$. Thus the subscript
$2$ in $E_2^{a,b}$ denotes a page, while the superscripts are degrees. In Cartan--Leray, $a$ is the group cohomology degree, $b$ is the
covering-space cohomology degree, and $a+b$ is the total degree in the quotient. The stable terms $E_\infty^{a,b}$ are the successive
quotients of a filtration of the cohomology group in that total degree; they do not generally give a canonical direct-sum decomposition.

For Cartan--Leray we use the construction in \cite[Chapter~VII, Section~7]{Br1982}. For a free $G$-cover $C\to Z$, let $EG$ be a contractible
free $G$-CW complex and $BG=EG/G$. The associated Borel fibration $C\to EG\times_G C\to BG$ has total space homotopy equivalent to $Z$. Its
cohomology spectral sequence has second page $H^a(G;H^b(C;A))$. This is the Cartan--Leray spectral sequence used below
\cite[Chapter~VII, Section~7, pp.~172--175]{Br1982}. For the cohomological indexing, convergence and edge maps, see also
\cite[Theorem~1.14]{HaSS}.

\begin{lemma}\label{lem6}
  For the cofinal regular covers in Proposition~\ref{prop1} and every abelian group $A$ with trivial $G$-action, there is a spectral sequence
  \begin{equation}\label{eq10}
    E_2^{a,b}
    =H^a\bigl(G;\Hinf^b(\R^n;A)\bigr)
    \ \Longrightarrow\ \Hinf^{a+b}(Y;A).
  \end{equation}
  Its only nonzero rows are $b=0,n-1$, both with trivial action. The edge maps are natural for coefficient homomorphisms and coefficient
  module cup products. In particular,
  \begin{equation}\label{eq11}
    H^q(G;A)\xrightarrow{\cong}\Hinf^q(Y;A)
    \quad(0\le q<n-1)
  \end{equation}
  and
  \begin{equation}\label{eq12}
    H^{n-1}(G;A)\longrightarrow\Hinf^{n-1}(Y;A)
    \quad\text{is injective}.
  \end{equation}
\end{lemma}

\begin{proof}
  Each free $G$-cover $C_j\to Z_j$ has a Cartan--Leray spectral sequence
  \[
    H^a(G;H^b(C_j;A))\Longrightarrow H^{a+b}(Z_j;A).
  \]
  One may obtain it from the Borel fibration, whose total space is homotopy equivalent to $Z_j$; see \cite[Chapter~VII, Section~7]{Br1982}. At
  each stage, $H^b(C_j;A)$ carries its induced $G$-module structure; no triviality of this action is assumed before taking the limit. The
  coverings and their $G$-actions are compatible under restriction.

  Take filtered direct limits. The bar resolution of a finite group has finite rank in each degree. Its cochain functors and cohomology
  therefore commute with filtered direct limits. Exactness of the limit also commutes with passage from one page of the spectral sequence to
  the next. For each fixed total degree the first quadrant filtration is finite, so the limiting spectral sequence converges to the direct
  limit of the abutments. Cofinality and \eqref{eq4} yield \eqref{eq10}. The upper row action is trivial by Proposition~\ref{prop1}.

  The first possible differential between the two rows is $d_n$. No incoming differential can reach a bottom-row position $a\le n-1$, since
  its source would have horizontal degree $a-n<0$. For total degrees below $n-1$ there is no upper-row term. This proves \eqref{eq11} and
  \eqref{eq12}. The bottom-row edge map at each stage is pullback by a classifying map $Z_j\to BG$ of the regular $G$-cover. These pullbacks
  are compatible under restriction and preserve coefficient maps and coefficient module cup products. Passing to the limit gives the claimed
  naturality of the edge maps.
\end{proof}

In rational coefficients, finite group cohomology vanishes in positive degrees. Thus \eqref{eq10} also gives
\begin{equation}\label{eq13}
  \Hinf^q(Y;\Q)=0\qquad(0<q<n-1).
\end{equation}
Only ordinary cohomology at infinity is used. No homotopy equivalence of an individual end neighborhood with a sphere quotient is claimed.

For the general transgression, see \cite[Chapter~VII]{Br1982} and \cite[Chapter~1, p.~54]{HaSS}. In this two-row spectral sequence the
\emph{transgression}  \cite{Br1982,HaSS} is $d_n:E_n^{0,n-1}\to E_n^{n,0}.$ All earlier differentials from these positions vanish, so its
source and target identify with $\Hinf^{n-1}(\R^n;\Z)$ and $H^n(G;\Z)$ in integral coefficients. It sends the chosen orientation generator of
the covering end to the obstruction class in degree $n$ on $BG$.

\begin{proposition}\label{prop2}
  If $G\ne1$, the transgression in \eqref{eq10} with integral coefficients gives a class
  \[
    e=d_n(1)\in H^n(G;\Z)
    \quad\text{with}\quad \ord(e)=|G|=d.
  \]
  Consequently $e$ is a unit in the Tate cohomology ring, and $n$ is a cohomological period of $G$.
\end{proposition}

\begin{proof}
  Lemma~\ref{lem4} and the two-row spectral sequence give the exact edge sequence
  \begin{equation}\label{eq14}
    0\longrightarrow H^{n-1}(G;\Z)
    \longrightarrow \Hinf^{n-1}(Y;\Z)
    \xrightarrow{f_\infty^*}\Hinf^{n-1}(\R^n;\Z)
    \xrightarrow{d_n}H^n(G;\Z)
    \longrightarrow\Hinf^n(Y;\Z)=0 .
  \end{equation}
  Both middle end groups are $\Z$. The connecting isomorphisms in \eqref{eq3} identify them with $H_c^n(Y;\Z)$ and $H_c^n(\R^n;\Z)$. Proper
  pullback is natural for this sequence, and has multiplication by $d$ on compact support top cohomology. Therefore $f_\infty^*$ is
  multiplication by $d$, up to the already fixed orientation convention. Since positive dimensional cohomology of a finite group is finite,
  the first group in \eqref{eq14} is zero. We obtain
  \[
    H^n(G;\Z)\cong\Z/d\Z,\qquad \ord(d_n(1))=d .
  \]

  The full order criterion for Tate periodicity \cite[Chapter~VI, Theorem~9.1]{Br1982} states that a positive degree integral class of order
  $|G|$ is invertible in Tate cohomology. Cup product with it gives periodicity for all coefficient modules. See also
  \cite[Section~1]{TW1971}. This criterion is algebraic and does not require a geometric realization of $G$ on a sphere.
\end{proof}

\subsection{Local periods}

Swan's $p$-period is described in \cite{Sw1960} and \cite[Section~1]{N1974}. We express it using units in the localized Tate cohomology ring,
with the Tate-periodicity conventions of \cite[Chapter~VI, Sections~5 and~9]{Br1982}.

\begin{definition}\cite{Sw1960,N1974,Br1982}\label{def4}
  Fix a prime $p$ dividing $|G|$, and let $R_p=\Z_{(p)}=\{a/b\in\Q:p\nmid b\}$. The \emph{minimal $p$-local period} $a_p$ is the least
  positive degree of a homogeneous unit in $\widehat H^*(G;R_p)$. A class $u$ of degree $a$ is such a unit if it has an inverse of degree $-a$
  whose product with $u$ is $1\in\widehat H^0(G;R_p)$. Equivalently, cup product with a periodicity unit in that degree gives isomorphisms
  $\widehat H^q(G;M)\to\widehat H^{q+a}(G;M)$ for every integer $q$ and every $R_pG$-module $M$.
\end{definition}

For a periodic finite group, such a degree exists by localizing an integral periodicity unit \cite{Br1982}. Multiplication adds the degrees of homogeneous
units, and inversion negates them. These degrees therefore form the subgroup $a_p\Z$ of $\Z$. Hence, if $n$ is a period, then
\begin{equation}\label{eq15}
  a_p\mid n.
\end{equation}
If $u\in H^{a_p}(G;R_p)$ is a unit, multiplication by $u$ identifies this group with
\[
  \widehat H^0(G;R_p)=R_p/|G|R_p .
\]
Writing $p^r$ for the order of a Sylow $p$-subgroup, we have
\begin{equation}\label{eq16}
  H^{a_p}(G;R_p)\cong\Z/p^r\Z.
\end{equation}
Here multiplication carries $1$ to $u$, and $R_p/|G|R_p\cong\Z/p^r\Z$ because the prime-to-$p$ factor of $|G|$ is invertible in $R_p$. In
particular, this group is indecomposable as an abelian group. These consequences of a Tate unit follow from
\cite[Chapter~VI, Sections~5 and~9]{Br1982}; the local-period convention is that of \cite{Sw1960,N1974}. The periods used here are minimal
local periods; an arbitrary degree of a unit need not divide another period.

\begin{lemma}\label{lem7}
  Every nontrivial finite group with periodic integral cohomology has a prime divisor $p$ for which $a_p=2$ or $a_p=4$.
\end{lemma}

\begin{proof}
  A finite periodic group has cyclic Sylow subgroups at odd primes, and its Sylow $2$-subgroup is cyclic or generalized quaternion. These are
  the standard periodicity criteria \cite[Chapter~VI]{Br1982}. If $|G|$ is even, Swan's formula gives $a_2=2$ in the cyclic case and $a_2=4$
  in the quaternion case \cite[Theorem~1]{Sw1960}; these formulas are also stated in \cite[Section~1, formula~(1.5)]{N1974}.

  If $|G|$ is odd, choose its smallest prime divisor $p$ and a cyclic Sylow subgroup $P$. Swan's other formula \cite[Theorem~2]{Sw1960} is
  \[
    a_p=2[N_G(P):C_G(P)].
  \]
  The quotient $N_G(P)/C_G(P)$ is a $p'$-subgroup of $\operatorname{Aut}(P)$: the Sylow subgroup $P$ already lies in the centralizer. Its
  order divides both $p-1$ and $|G|$. Every prime dividing $p-1$ is smaller than $p$, so minimality of $p$ makes this quotient trivial.
  Consequently $a_p=2$.
\end{proof}

The quaternion formula includes fusion in the full group, rather than only an assertion about a quaternion subgroup. This is why the local
period theorem, not merely the classification of Sylow subgroups, is used in Lemma~\ref{lem7}.

\section{An ambient Bockstein obstruction}\label{sec5}

\subsection{Classes extending over the two sides}

The complement $D$ in \eqref{eq5} may be wild. We therefore use neighborhoods of $D$ and the continuity property of \v{C}ech cohomology
\cite{Sp1966}. For Mayer--Vietoris we use the cohomological convention in \cite[Sections~2.2 and~3.1]{Ha2002}. For an open cover $X=U\cup V$,
the cohomological \emph{Mayer--Vietoris sequence} \cite{Ha2002} is the natural long exact sequence 
\[
  \cdots\to H^q(X;A)\to H^q(U;A)\oplus H^q(V;A)
  \to H^q(U\cap V;A)\to H^{q+1}(X;A)\to\cdots .
\]
The middle map is the difference of the two restrictions. We use this sequence for open neighborhoods even when one limiting side is a wild
compact set.

\begin{lemma}\label{lem8}
  For every abelian coefficient group $A$ there are natural isomorphisms
  \begin{equation}\label{eq17}
    H^q(Y;A)\oplus\Cech^q(D;A)
    \xrightarrow{\cong}\Hinf^q(Y;A),
    \qquad 1\le q\le n-2.
  \end{equation}
  The maps of each summand are restrictions. They preserve cup products, and more generally coefficient module cup products, whenever the
  products are defined. Their degree-$(n-1)$ source groups vanish.
\end{lemma}

\begin{proof}
  For a compact set $K\subset Y$, put $U_K=\Sph^n\setminus K$. The open cover $\Sph^n=Y\cup U_K$ has intersection $Y\setminus K$. For
  $1\le q\le n-2$, the neighboring cohomology groups of $\Sph^n$ in Mayer--Vietoris vanish, so restriction gives
  \[
    H^q(Y;A)\oplus H^q(U_K;A)
    \xrightarrow{\cong}H^q(Y\setminus K;A).
  \]
  We choose the plus sign for the two restrictions; changing the usual Mayer--Vietoris sign does not alter this isomorphism.

  The $U_K$ form a cofinal system of open neighborhoods of $D$. On these manifold neighborhoods singular and \v{C}ech cohomology agree.
  Continuity gives $\varinjlim_K H^q(U_K;A)=\Cech^q(D;A)$. Taking limits proves \eqref{eq17}. Both summand maps come from restrictions of
  cochains, so they commute with coefficient homomorphisms and the indicated products. Their definitions persist in degree $n-1$, even though
  the direct-sum isomorphism does not. The source groups in that degree vanish by Lemma~\ref{lem4}.
\end{proof}

\begin{remark}\label{rem2}
  Equation~\eqref{eq17} is an additive decomposition. We do not claim that mixed products of the two summands vanish or that it is an
  isomorphism with a direct-product ring. The argument uses only that classes extending over a fixed side retain that property under
  coefficient operations and products.
\end{remark}

For the description of the Pr\"ufer group we use \cite{Fu1970}. Fix $p\mid d$, and for the remainder of this section put
\begin{equation}\label{eq18}
  R=\Z_{(p)},\qquad T=\Q/R .
\end{equation}
The \emph{Pr\"ufer $p$-group} $T$ is the union of its cyclic subgroups generated by $1/p^m+R$, of order $p^m$. Equivalently, it is the direct
limit of $\Z/p^m\Z$ under the maps sending a residue class to $p$ times that class in $\Z/p^{m+1}\Z$. It is a divisible $p$-primary torsion
group \cite{Fu1970}. Here it is considered as an $R$-module, not as a coefficient ring. We use the Bockstein convention in
\cite[Section~3.E]{Ha2002} and \cite[Chapter~III]{Br1982}. The exact sequence
\begin{equation}\label{eq19}
  0\longrightarrow R\longrightarrow\Q\longrightarrow T
  \longrightarrow0
\end{equation}
defines Bockstein homomorphisms $\beta:H^q(-;T)\to H^{q+1}(-;R)$, including their \v{C}ech, group cohomology and end-cohomology forms. They
are the connecting homomorphisms of the long exact coefficient sequence \cite{Ha2002,Br1982}. At the cochain level, lift a $T$-valued cocycle
$c$ to a $\Q$-valued cochain $\widetilde c$. Its coboundary is $R$-valued, and its class is $\beta[c]$. We call a class $v$ with $\beta(v)=u$
a \emph{Bockstein preimage} of $u$. Such a class exists when $u$ lies in the image of $\beta$ and is unique when $\beta$ is injective.
Lemma~\ref{lem9} specifies the degrees in which the Bocksteins used here are isomorphisms.

If $x$ has $T$ coefficients and $y$ has $R$ coefficients, the module cup product satisfies
\[
  \beta(x\smile y)=\beta(x)\smile y.
\]
Indeed, an $R$-valued cocycle representing $y$ has zero coboundary, so the cochain Leibniz rule gives
$\delta(\widetilde c\smile y)=\delta\widetilde c\smile y$. This is the Bockstein product formula in the coefficient form needed below; see
\cite[Sections~3.2 and~3.E]{Ha2002} and \cite[Chapter~V, Section~3]{Br1982}. The same identity passes to the neighborhood limits defining
\v{C}ech and end cohomology.

\begin{lemma}\label{lem9}
  The Bockstein in \eqref{eq19} is an isomorphism
  \[
    H^q(Y;T)\cong H^{q+1}(Y;R),\qquad
    \Cech^q(D;T)\cong\Cech^{q+1}(D;R)
    \quad(q>0).
  \]
  It is also an isomorphism
  \[
    \Hinf^q(Y;T)\cong\Hinf^{q+1}(Y;R)
    \quad(1\le q\le n-3).
  \]
  These maps commute with the restrictions in \eqref{eq17} and the Cartan--Leray edge maps.
\end{lemma}

\begin{proof}
  For the two sides, rational cohomology in positive degrees vanishes by Lemma~\ref{lem4}. For the end, use \eqref{eq13}. The long exact
  coefficient sequences give the asserted isomorphisms. Coefficient sequences in \v{C}ech cohomology may be obtained from the neighborhood
  limits; filtered direct limits preserve their exactness. Naturality gives compatibility with restrictions and edge maps.
\end{proof}

\subsection{Excluding a smaller local period}

\begin{proposition}\label{prop3}
  Suppose $G\ne1$, and let $a=a_p$ be a minimal local period. If $a\mid n$ and $2\le a\le n-2$, the preceding end covering cannot exist. In
  particular, for an even $n$, every minimal local period of the end group must equal $n$.
\end{proposition}

\begin{proof}
  Choose a periodicity unit $u\in H^a(G;R)$ and use \eqref{eq11} to regard it as $u_\infty\in\Hinf^a(Y;R)$. The entire group $\Hinf^a(Y;R)$ is
  cyclic of order $p^r$ by \eqref{eq16}. Its additive decomposition \eqref{eq17} is a direct sum. A finite cyclic group of prime-power order
  has no nontrivial direct-sum decomposition. Thus one summand is zero, and $u_\infty$ comes wholly from the other side. Denote that side by
  $S$, where cohomology means ordinary cohomology if $S=Y$ and \v{C}ech cohomology if $S=D$. There is a class $u_S\in H^a(S;R)$ whose
  restriction is $u_\infty$.

  Finite-group cohomology with rational coefficients vanishes in positive degrees. The group Bockstein therefore gives a unique class
  $\widehat u\in H^{a-1}(G;T)$ with $\beta_G(\widehat u)=u$. Since $a\ge2$, Lemma~\ref{lem9} also gives a class $\widehat u_S\in H^{a-1}(S;T)$
  with $\beta_S(\widehat u_S)=u_S$. The end Bockstein is an isomorphism in these degrees. Naturality consequently identifies the restriction
  of $\widehat u_S$ with the edge image $\widehat u_\infty$ of $\widehat u$.

  Write $k=n/a$. Form the coefficient module cup product
  \begin{equation}\label{eq20}
    \alpha=\widehat u\smile u^{k-1}
    \in H^{n-1}(G;T).
  \end{equation}
  The coefficient sequence \eqref{eq19} is a sequence of $R$-modules. Since the second factor has $R$ coefficients, the Bockstein product
  formula \cite[Sections~3.2 and~3.E]{Ha2002} gives
  \begin{equation}\label{eq21}
    \beta_G(\alpha)=u^k .
  \end{equation}
  The power $u^k$ is a Tate unit, hence is nonzero in the positive degree group $H^n(G;R)$. Therefore $\alpha\ne0$. By the injective edge
  \eqref{eq12}, its image $\alpha_\infty\in\Hinf^{n-1}(Y;T)$ is nonzero.

  On the other hand, compatibility of the edge with module cup products gives
  \[
    \alpha_\infty
    =\widehat u_\infty\smile u_\infty^{k-1}
    =\operatorname{res}_S
    \bigl(\widehat u_S\smile u_S^{k-1}\bigr).
  \]
  This product extends over the single side $S$. But $H^{n-1}(Y;T)=\Cech^{n-1}(D;T)=0$ by \eqref{eq7}. It follows that $\alpha_\infty=0$, a
  contradiction.

  For a periodic finite group each minimal local period is even \cite{Sw1960,Br1982}. If $n$ is even and $a<n$ with $a\mid n$, then
  $a\le n-2$. Proposition~\ref{prop2} and \eqref{eq15} prove the last assertion.
\end{proof}

\begin{remark}\label{rem3}
  The nonvanishing in this proof is established in group cohomology before passing through the injective edge map. The end Bockstein in degree
  $n-1$ cannot be used to prove it: $\Hinf^n(Y;R)=0$. Similarly, the product in \eqref{eq20} uses the $R$-module structure on $T$, not a ring
  structure on $T$.
\end{remark}

\section{Higher-dimensional injectivity}\label{sec6}

\begin{proof}[Proof of Theorem~\ref{thm1}]
  Replace $B^n$ by $\R^n$ using a source homeomorphism. All the hypotheses, including compactness of the branch set, are preserved, because a
  source homeomorphism preserves local homeomorphisms and compact inverse images. Fix the orientation so that $d>0$. Lemma~\ref{lem5} provides
  connected end neighborhoods whose full inverse images are $d$-sheeted coverings. Proposition~\ref{prop1} makes a cofinal sequence of these
  coverings regular, with compatible deck group $G$ and $|G|=d$.

  Suppose that $d>1$. The group $G$ is then nontrivial. By Lemma~\ref{lem6}, the cohomology at infinity is computed by the two-row
  Cartan--Leray spectral sequence. Its integral edge sequence \eqref{eq14}, together with proper degree theory, gives a class $e\in H^n(G;\Z)$
  of order $d=|G|$. Proposition~\ref{prop2} makes $e$ a Tate unit, so $n$ is a cohomological period of $G$. After localization at each prime
  dividing $d$, the minimal period divides $n$ by \eqref{eq15}.

  Lemma~\ref{lem7} supplies a prime $p\mid d$ whose minimal local period $a=a_p$ is $2$ or $4$. If $n$ is odd, the divisibility $a\mid n$ is
  impossible. If $n$ is even, then $n\ge6$ and $2\le a\le n-2$. Proposition~\ref{prop3} excludes this case: the cyclic $p$-primary periodicity
  class extends over one side of the ambient decomposition \eqref{eq17}; its Bockstein preimage extends over that side as well. Their product
  gives a nonzero group class in degree $n-1$ and hence a nonzero end class by \eqref{eq12}. But the same product extends over a side whose
  degree-$(n-1)$ cohomology vanishes by \eqref{eq7}. This is the required contradiction.

  Consequently $d=1$. For each $y\in Y$, \eqref{eq2} expresses $1$ as the sum of positive local indices over the nonempty fiber, so that fiber
  has exactly one point. The map is therefore a continuous open bijection onto $Y$, and hence a homeomorphism.
\end{proof}
\begin{remark}
The argument also recovers the three-dimensional conclusion of \cite{KLT2022}: a nontrivial end group would have a minimal local period two or
four dividing three. The three-dimensional result is credited to that paper; the ambient Bockstein obstruction is needed for the even higher
dimensions.
\end{remark}
\begin{corollary}\label{cor1}
  Let $n\ge5$ and let $f:B^n\to f(B^n)\subset\R^n$ be a proper nonconstant quasiregular map. If $B_f$ is compact, then $f$ is a quasiconformal
  homeomorphism onto its image.
\end{corollary}

\begin{proof}
  By Reshetnyak's theorem, we get that $f$ is openness and discreteness. Hence,  Theorem~\ref{thm1} gives the homeomorphism conclusion. A homeomorphic quasiregular
  map is quasiconformal.
\end{proof}

\begin{corollary}\label{cor2}
  Let $n\ge5$, let $\Omega\subset\R^n$ be a domain, and let $f:\Omega\to\R^n$ be a branched covering. Suppose that $B(a,r)\Subset\Omega$ is a
  normal domain and $B_f\cap\partial B(a,r)=\varnothing$. Then $f|_{B(a,r)}$ is injective.
\end{corollary}

\begin{proof}
  The normal domain identity makes the restriction proper onto $f(B(a,r))$. Since $B_f$ is closed in $\Omega$ and misses the boundary sphere,
  its intersection with the ball is compact in the ball. Theorem~\ref{thm1} applies after a translation and scaling.
\end{proof}

Corollary~\ref{cor2} gives the higher-dimensional normal ball version of \cite[p.~125, Open problem~9.18]{V1988}, whose three-dimensional case
was deduced in \cite{KLT2022}. It also rules out introducing a nonempty branch set inside such a normal ball while retaining a branch free
boundary.

\begin{remark}\label{rem4}
  Both properness and the compact branch hypothesis are necessary. In cylindrical coordinates the winding map
  \[
    (r,\theta,x_3,\ldots,x_n)
    \longmapsto(r,m\theta,x_3,\ldots,x_n),\qquad m\ge2,
  \]
  is proper and noninjective, but its branch set is an unbounded $(n-2)$-plane. Conversely, an exponential map on a bounded rectangular strip,
  multiplied by identity maps, is a noninjective local homeomorphism on a domain homeomorphic to a ball and fails properness onto its image;
  see \cite{KLT2022}. The planar map $z\mapsto z^m$ has compact branch set and is proper, so dimension two cannot be included in
  Theorem~\ref{thm1}.
\end{remark}

\section{The four-dimensional quasiregular case}\label{sec7}

The argument in this section uses the analytic hypothesis in an essential way. Throughout the section, let
\[
  f:B^4\longrightarrow Y\subset\R^4
\]
be a proper nonconstant quasiregular mapping whose branch set is compact, and let $d=\deg f$. Quasiregularity is understood in the sense
specified above; in particular, the distortion bound is uniform on $B^4$. The topological consequences of quasiregularity give an oriented
branched covering of positive degree, and every restriction of $f$ that is a homeomorphism is quasiconformal \cite{R1993,IM2001}.

\subsection{Quasiconformal structures and four-dimensional topology}

The metric-space definition of quasisymmetry is given in \cite[Section~1.1]{TV1980}.
For the Euclidean quasiconformal background, see \cite{Va1971}.

\begin{definition}\cite{TV1980,Va1971}\label{def5}
  A homeomorphism $h:(X,d_X)\to (Z,d_Z)$ is \emph{quasisymmetric} if there is an increasing homeomorphism
  $\eta:[0,\infty)\to[0,\infty)$ such that
  \[
    \frac{d_Z(h(x),h(a))}{d_Z(h(x),h(b))}
    \le \eta\!\left(\frac{d_X(x,a)}{d_X(x,b)}\right)
  \]
  for every triple of distinct points $x,a,b\in X$.
\end{definition}

For locally quasiconformal manifolds with boundary we follow \cite[Section~2.1]{TV1981}.
The atlas definition for closed four-manifolds is given in \cite[Definition~4.1(i)]{DS1989}.

\begin{definition}\cite{TV1981,DS1989}\label{def6}
  A \emph{locally quasiconformal manifold}, or \emph{LQC manifold}, is a topological manifold
  with an atlas whose coordinate changes are locally quasiconformal.

  For a manifold with boundary the charts take values in relatively open subsets of a closed Euclidean half-space.
  At every point of a transition domain, including a boundary point, there is a relative neighborhood
  on whose interior the transition is quasiconformal with a finite uniform bound.
  The transition extends continuously to the boundary. Boundary restrictions give the induced LQC structure in one lower dimension.
\end{definition}

Brown defines collars and local flatness in \cite[Sections~II and~IV]{Br1962}.

\begin{definition}\cite{Br1962}\label{def7}
  An embedding of a codimension-one manifold is \emph{locally flat} if each point has an ambient topological coordinate neighborhood in which
  the embedded manifold is a coordinate hyperplane.

  A \emph{topological collar} of $\partial W$ is an embedding $\partial W\times[0,\varepsilon)\to W$ restricting to the identity on
  $\partial W\times\{0\}$, whose image is a neighborhood of the boundary. These topological notions do not themselves specify a quasiconformal
  structure.
\end{definition}

The analytic boundary extension used below is \cite[Theorem~4.10]{TV1981}. For cohomology manifolds we follow \cite[Chapter~I, Definition~3.3]{Bo1960}.
Their role in Smith theory is discussed in \cite[Chapter~V]{Bo1960} and \cite{Bd1972}.

\begin{definition}\cite{Bo1960,Bd1972}\label{def8}
  A \emph{cohomology manifold over a coefficient ring} is a finite-dimensional locally compact space whose local cohomology and orientation
  sheaf have the properties of those of Euclidean space over that ring. The coefficient rings used below are $\Z$ and $\mathbb F_p$.
  Topological manifolds without boundary are cohomology manifolds over either ring. Local fixed sets in Smith theory are understood in this
  sense; no local linearity or manifold structure over $\Z$ is asserted for such fixed sets.
\end{definition}

We also need the following conceptions. 
\begin{definition}\label{def9}
  The \emph{binary icosahedral group} $I^*$ is the inverse image of the rotation group of the icosahedron under
  $\mathrm{SU}(2)\to\mathrm{SO}(3)$. It has order $120$ and acts freely by left multiplication on the unit quaternions
  $S^3$ \cite[p.~75]{Ha2002}. The group $I^*$ is isomorphic to $\mathrm{SL}_2(\mathbb F_5)$ \cite[p.~222]{Hi2002}.
  The quotient $S^3/I^*$ is the \emph{Poincar\'e homology sphere} \cite[Supplement~10, p.~855]{C1978}.

  A \emph{spherical space form} is a quotient of a round sphere by a finite group of isometries acting freely
  \cite[Sections~7.1 and~11.5]{Hi2002}.
  An \emph{integral homology ball} is a compact connected manifold with nonempty boundary and the integral homology of a point
  \cite[Section~76.2, p.~1771]{FriedlAT}.
  An \emph{integral homology sphere} is a closed manifold with the integral homology of a sphere of the same dimension
  \cite[p.~833]{C1978}.
\end{definition}

The closed-embedding characterization of absolute neighborhood retracts is given in \cite[Theorem~12.1]{Dug1958}.

\begin{definition}\cite[Theorem~12.1]{Dug1958}\label{def10}
  A metrizable space $X$ is an \emph{absolute neighborhood retract} (ANR) if, whenever $X$ is realized as a closed subset of a metrizable
  space $Z$, some open neighborhood of $X$ in $Z$ retracts onto $X$.

  Manifolds, including manifolds with boundary, are ANRs. Open second countable manifolds have CW homotopy type \cite{Mi1959}. For a compact
  ANR, West's theorem \cite[Corollary~5.3]{We1977} gives the homotopy type of a finite CW complex.
\end{definition}

We use the cohomological definition of the intersection form in \cite[Definition~1.2.1]{GS1999} and the lattice terminology of \cite{MH1973}.

\begin{definition}\cite{GS1999,MH1973}\label{def11}
  The \emph{intersection form} of a closed oriented four-manifold $X$ is the symmetric integral bilinear form
  \[
    Q_X(a,b)=\langle a\smile b,[X]\rangle
    \quad\text{on } H^2(X;\Z)/\operatorname{Tor} H^2(X;\Z).
  \]
  It is \emph{even} if $Q_X(a,a)$ is even for every $a$, and \emph{definite} if its extension to real coefficients is positive or negative
  definite. A lattice is \emph{unimodular} if the determinant of a Gram matrix is $\pm1$.
\end{definition}

Plumbing is defined in \cite[Example~4.6.2]{GS1999}.
The negative $E_8$ plumbing and its boundary are described in \cite[Figure~4.33 and Exercise~5.1.12(a)]{GS1999}.
For the $E_8$ lattice, see \cite{MH1973}.

\begin{definition}\cite{GS1999,MH1973}\label{def12}
  \emph{Plumbing} oriented disk bundles over oriented surfaces means identifying small products $D^2\times D^2$ in two bundles by exchanging
  the base and fiber factors.

  The $E_8$ plumbing uses eight bundles over $S^2$ with Euler number $-2$, with one plumbing for each edge of the $E_8$ Dynkin tree. This tree
  has eight vertices and seven edges, with one trivalent vertex and three arms of lengths $1$, $2$ and $4$ edges. A vertex denotes a disk
  bundle, and an edge denotes a plumbing between the corresponding bundles. With the standard orientations its intersection matrix has
  diagonal entries $-2$ and entry $1$ for adjacent vertices; all other entries are zero. Its intersection form is $-E_8$ and its boundary is a
  Poincar\'e homology sphere, with the orientation determined by the plumbing.
\end{definition}

\subsection{Compactifying the quasiconformal end}

\begin{proposition}\label{prop4}
  There is a compact connected oriented LQC four-manifold $W$ with
  \[
    \operatorname{int}W=Y,\qquad \partial W=S^3/G,
  \]
  where $G$ is the regular end-covering group in Proposition~\ref{prop1}. The extended action of $G$ on $S^3$ is free and orientation
  preserving. Consequently, $S^3\to\partial W$ is a regular covering of degree $d$.
\end{proposition}

\begin{proof}
  Apply Proposition~\ref{prop1} to obtain a connected neighborhood $V$ of infinity in $Y$ such that
  \[
    C=f^{-1}(V)\longrightarrow V
  \]
  is a regular $d$-sheeted covering with deck group $G$. The set $C$ contains an entire outer annulus of $B^4$ and is disjoint from a closed
  ball containing $f^{-1}(f(B_f))$. Its Euclidean boundary is therefore the disjoint union of $S^3$ and a compact subset of $B^4$. Every deck
  transformation $g:C\to C$ is quasiconformal. Indeed, near each point it is the composition of a local inverse branch of $f$ with another
  local branch, and the resulting distortion bound depends only on the distortion of $f$ and the dimension \cite{R1993,IM2001}. The deck
  transformations are orientation preserving, since $f\circ g=f$ and the unbranched local branches of $f$ preserve orientation.

  We first verify that every deck transformation preserves the outer boundary. If $x_j\in C$ and $|x_j|\to1$, properness implies that $f(x_j)$
  eventually leaves every compact subset of $Y$. If a subsequence of $g(x_j)$ converged to a point of $B^4$, continuity of $f$ at that point
  would contradict $f(g(x_j))=f(x_j)$. Thus $|g(x_j)|\to1$. Applying the same argument to $g^{-1}$ shows that the cluster set of $g$ at $S^3$
  is exactly $S^3$. We prove the required spherical boundary extension directly. Let $R(x)=x/|x|^2$ be reflection in $S^3$ and put
  \[
    C^*=C\cup S^3\cup R(C).
  \]
  This is an open neighborhood of $S^3$ in $\R^4$. Choose $r<1$ so that the annulus $A=\{x:r<|x|<1\}$ lies in $C$ and $|g(x)|>1/2$ on $A$. The
  latter follows from the outer-end preservation just proved. Let $K_g<\infty$ be an outer dilatation bound for $g$. Quasiconformality and the
  change-of-variables formula give
  \[
    \int_A\|Dg\|^4\,dx
    \le K_g\int_A J_g\,dx
    =K_g|g(A)|\le K_g|B^4|<\infty.
  \]
  Since $g$ is bounded, it belongs to $W^{1,4}(A)$ and has a Sobolev trace $\tau_g$ on $S^3$. Outer-end preservation is uniform as $|x|\to1$,
  so $|\tau_g|=1$ almost everywhere on $S^3$. Both inversions in $R\circ g\circ R$ have bounded derivatives on the relevant annuli. Thus this
  map belongs to $W^{1,4}(R(A))$, and its trace is $R\circ\tau_g=\tau_g$. Equality of the traces permits Sobolev gluing across the smooth
  sphere, because the interface terms cancel in the weak derivative identity. The pasted map belongs to $W^{1,4}$ near $S^3$ and satisfies the
  bounded-distortion inequality almost everywhere: the inequality holds on both open sides, and $S^3$ has four-dimensional measure zero.
  Reshetnyak's theorem \cite{Re1989,R1993} gives a continuous representative agreeing with the prescribed maps on those sides. Outer-end
  preservation then makes its boundary image lie in $S^3$.

  Apply the same construction to $g^{-1}$. The two continuous extensions take $C^*$ into itself and are inverse on the dense subset
  $C^*\setminus S^3$. Their compositions are therefore the identity everywhere. Consequently the extension $\widehat g:C^*\to C^*$ is a
  homeomorphism, is quasiconformal across $S^3$, and restricts to a homeomorphism of $S^3$. It is given explicitly by
  \[
    \widehat g(x)=
    \begin{cases}
      g(x),&x\in C,\\
      g(x),&x\in S^3,\\
      R\circ g\circ R(x),&x\in R(C),
    \end{cases}
  \]
  where the middle expression denotes the boundary extension. The explicit formula and continuity on $S^3$ show that these extensions satisfy
  \[
    \widehat g\circ\widehat h=\widehat{gh}
    \quad\text{on }C^*.
  \]
  They preserve orientation. The action is free on $C$ and $R(C)$.

  We claim that it is also free on $S^3$. Suppose a boundary point had a nontrivial stabilizer. By Cauchy's theorem, the stabilizer contains
  an element $h$ of prime order $p$. The local fixed set of $\widehat h$ in $C^*$ equals its local fixed set in $S^3$, because the two sides
  have no fixed points. Local Smith theory makes a nonempty local fixed set a mod-$p$ cohomology manifold. For odd $p$, its dimension has the
  same parity as the ambient dimension. For $p=2$, the same parity conclusion holds for an orientation-preserving involution on an orientable
  integral cohomology manifold \cite[Chapter~V, Theorems~2.2, 2.3 and~2.5]{Bo1960}; see also \cite[proof of Lemma~2.3]{Ye2018} for the
  involution form used here. In the four-dimensional ambient manifold, the fixed-set germ therefore has even dimension. If $h|_{S^3}$ were the
  identity, this germ would be the three-dimensional germ of $S^3$, since there are no fixed points on either open side. This is already a
  contradiction. Otherwise $h|_{S^3}$ has order $p$, so the same Smith theorems apply to this effective boundary action. The boundary action
  preserves the induced orientation of $S^3$, because the ambient action preserves orientation and the two sides of the sphere separately. The
  common local fixed set would consequently have odd dimension when regarded in $S^3$ and even dimension when regarded in $C^*$. These are
  dimensions of the same local mod-$p$ cohomology manifold germ, an intrinsic invariant of that germ. This is impossible.

  We next construct the compactification without making an assumption on the Euclidean boundary of $Y$. On $\overline B^4$, identify two
  interior points precisely when they have the same image under $f$, and identify two boundary points precisely when they lie in the same
  $G$-orbit. No interior point is identified with a boundary point. This equivalence relation is closed. For limits of pairs of interior
  points, continuity settles the case in which both limits remain interior. A pair cannot have one interior and one boundary limit, by
  properness of $f$. If both limits belong to $S^3$, the corresponding points eventually lie in $C$. Equal images there differ by a deck
  transformation; since $G$ is finite, a subsequence uses one fixed deck transformation, which relates the two boundary limits. Limits of
  boundary pairs are handled directly by finiteness and continuity of the action. The quotient $W$ is therefore compact and Hausdorff. Its
  interior is $Y$ with its original topology, since $f$ is an open quotient map. The subset $C\cup S^3$ of $\overline B^4$ is open and
  saturated, because $C=f^{-1}(V)$. Its quotient is a neighborhood of the new boundary.

  To describe its boundary charts, choose a small neighborhood of a point of $S^3$ in $C^*$ whose translates under distinct elements of $G$
  are disjoint. Such neighborhoods exist because the boundary action is free. Their quotients, restricted to the inside of the sphere, are
  topological half-space charts. Smooth sphere-flattening coordinates make transitions between these charts restrictions of smooth conjugates
  of the quasiconformal mappings $\widehat g$. On the flattened boundary these restrictions are locally quasisymmetric: the ambient
  quasiconformal maps are locally quasisymmetric, and this inequality restricts to boundary triples. They are therefore locally quasiconformal
  in dimension three \cite{Va1971,TV1980}. Their overlaps with interior Euclidean charts on $Y$ have coordinate changes given by the local
  quasiconformal branches of $f$ or their inverses. Thus they form an LQC atlas with boundary. The orientation is compatible with the
  orientation of $Y$, and $\partial W=S^3/G$ is connected. The free boundary action gives the asserted regular covering.
\end{proof}

\subsection{Homology of the compact filling}

\begin{lemma}\label{lem10}
  The compact manifold $W$ in Proposition~\ref{prop4} is simply connected and has the integral homology of a point. Its boundary is an
  integral homology three-sphere.
\end{lemma}

\begin{proof}
  An ordinary topological collar of the boundary shows that $Y=\operatorname{int}W\hookrightarrow W$ is a homotopy equivalence
  \cite[Theorem~2]{Br1962}. Lemma~\ref{lem1}, applied to the proper mapping $f$, gives
  \[
    \widetilde H_i(Y;\Q)=0.
  \]
  Hence $W$ is a rational homology ball.

  Lemma~\ref{lem2} applies after a topological identification of $B^4$ with $\R^4$, and proves that $\pi_1(Y)\cong\pi_1(W)$ is finite. This
  identification is used only for the covering space argument, not to transport the quasiregular structure. Put $k=|\pi_1(W)|$ and let
  $\widetilde W\to W$ be the universal covering. It is a finite covering of a compact topological manifold, so $\widetilde W$ is compact. Its
  interior is the universal covering $N$ of $Y$. Lemma~\ref{lem2} also supplies a proper surjective branched covering of $N$ by $B^4$, of some
  positive degree. A second application of Lemma~\ref{lem1} gives
  \[
    \widetilde H_i(N;\Q)=0.
  \]
  The inclusion $N\hookrightarrow\widetilde W$ is a homotopy equivalence by a boundary collar. Compact topological manifolds are compact ANRs,
  and compact ANRs have the homotopy type of finite CW complexes \cite[Corollary~5.3]{We1977}. Consequently their rational homology determines
  the Euler characteristics, and
  \[
    \chi(W)=\chi(\widetilde W)=1.
  \]
  Euler characteristic is therefore multiplicative under this finite covering \cite[Section~2.2, Exercise~22]{Ha2002}. To apply the finite-CW
  formula, pull the covering back along a finite CW model homotopy equivalent to $W$. The pulled-back covering is homotopy equivalent to
  $\widetilde W$, and counting the lifted cells gives the formula. Thus no triangulation of $W$ is assumed:

  \[
    1=\chi(\widetilde W)=k\chi(W)=k.
  \]
  Thus $W$ is simply connected.

  Set $M=\partial W$. The finite oriented covering $S^3\to M$ and its transfer make $M$ a rational homology three-sphere
  \cite[Section~3.G]{Ha2002}. In particular, $H_1(M;\Z)$ is finite and, by three-dimensional duality, $H_2(M;\Z)=0$. All positive dimensional
  integral homology groups of $W$ are finitely generated torsion groups, since $W$ is a rational homology ball. Simple connectivity gives
  $H_1(W;\Z)=0$. The universal coefficient theorem then yields
  \[
    H^1(W;\Z)=H^2(W;\Z)=0.
  \]
  Poincar\'e--Lefschetz duality  \cite[Section~3.3]{Ha2002}, valid for oriented topological manifolds, gives
  \[
    H_3(W,M;\Z)=0,\qquad H_2(W,M;\Z)=0
  \]
 In the homology sequence of $(W,M)$, the map $H_4(W,M;\Z)\to H_3(M;\Z)$ sends the relative fundamental class to
  the boundary fundamental class and is an isomorphism, because $M$ is connected. It follows that $H_3(W;\Z)=0$. The remaining sequence, using
  $H_2(M;\Z)=0$, gives
  \[
    H_2(W;\Z)=H_1(M;\Z)=0.
  \]
  Finally $H_4(W;\Z)=0$ since $\partial W$ is nonempty. These conclusions prove the assertions.
\end{proof}

\begin{lemma}\label{lem11}
  If $d>1$, then $G\cong I^*=\operatorname{SL}(2,\mathbb F_5)$ and $M=\partial W$ is homeomorphic to the Poincar\'e homology three-sphere,
  with one of its two orientations.
\end{lemma}

\begin{proof}
  The covering $S^3\to M$ in Proposition~\ref{prop4} is universal, and hence $\pi_1(M)\cong G$. By Lemma~\ref{lem10}, $H_1(M;\Z)=0$, so $G$ is
  a finite perfect group. The binary icosahedral group $I^*$ is the only nontrivial finite group that can be the fundamental group of an
  integral homology three-sphere; equivalently, it is the only nontrivial finite perfect group of cohomological period four
  \cite[Section~14.13, p.~292]{Hi2002}. Thus $G\cong I^*$. The classification of spherical three-manifolds identifies the resulting integral
  homology three-sphere as the Poincar\'e sphere $P=S^3/I^*$, up to orientation. Here one may use elliptization together with the
  classification of spherical space forms; the relevant finite fundamental-group consequence of geometrization is recorded in
  \cite[Remark~1.4]{Pe2003}.
\end{proof}

\subsection{The quasiconformal gluing obstruction}

\begin{lemma}\label{lem12}
  There is no compact simply connected LQC four-manifold with the integral homology of a point and boundary homeomorphic to the Poincar\'e
  homology three-sphere.
\end{lemma}

\begin{proof}
  Suppose that $W$ were such a manifold and write $M=\partial W$. Equip the standard Poincar\'e sphere $P$ with its smooth LQC structure. The
  boundary $M$ has the induced three-dimensional LQC structure. By the three-dimensional LQC Hauptvermutung, a homeomorphism $M\to P$ can be
  approximated by an LQC homeomorphism \cite[Corollary~4.5]{TV1981}
  \[
    h:M\longrightarrow P.
  \]
  Choose the approximation close enough to preserve the required orientation.

  Let $E$ be the smooth four-manifold obtained by plumbing disk bundles with Euler number $-2$ along the $E_8$ diagram. It is simply
  connected, has intersection form $-E_8$, and its boundary is the Poincar\'e sphere, up to orientation \cite{GS1999}. If necessary reverse
  the orientation of $E$ so that $h$ can be used for an oriented gluing. The smooth boundary collar of $E$ supplies a standard smooth product
  $P\times[0,\varepsilon)$. Apply the boundary-extension theorem \cite[Theorem~4.10]{TV1981} to extend $h$ to an LQC embedding of a
  neighborhood of $M$ in $W$ into this product. This theorem explicitly includes dimension four. Its extension preserves the boundary and
  takes interior points to interior points, as in the clean embeddings used in \cite[Theorems~4.8 and~4.10]{TV1981}. Doubling the extension
  across the boundary gives an injective map between four-manifolds without boundary. Invariance of domain makes its image open; restricting
  back to the half-space shows that the original image is a relatively open neighborhood of $P\times\{0\}$. Since $P$ is compact, it contains
  $P\times[0,\delta)$ for some $\delta>0$. Replace the domain by the inverse image of this product and relabel $\delta$ as $\varepsilon$. The
  extension is then an LQC homeomorphism
  \[
    H:U\longrightarrow P\times[0,\varepsilon).
  \]

  Form the closed oriented topological manifold
  \[
    X=W\cup_h E.
  \]
  The map $H$ and the smooth collar on the $E$ side give coordinate neighborhoods across the seam modeled on
  $P\times(-\varepsilon,\varepsilon)$. In local smooth coordinates on $P$, transitions between two such seam charts are smooth product
  coordinate changes. Transitions with an interior chart on the $W$ side are locally quasiconformal by the defining property of $H$;
  transitions with interior charts on the $E$ side are smooth. Thus $X$ has a quasiconformal atlas. Compactness permits a finite cover by
  charts $V_i$ whose closures are compact subsets of larger chart domains $O_i$. The closure of $V_i\cap V_j$ lies in $O_i\cap O_j$. A finite
  cover of this compact set by neighborhoods with local distortion bounds gives one finite bound for the transition on each restricted
  overlap. The resulting compatible maximal quasiconformal atlas satisfies \cite[Definition~4.1]{DS1989}. The product used here is the fixed
  smooth product on $P$; no product construction for arbitrary LQC spaces is being assumed. The actual four-dimensional map $H$, furnished by
  the boundary extension theorem \cite[Theorem~4.10]{TV1981}, provides the needed compatibility.

  Take open neighborhoods of the two pieces whose intersection is a collar of $M$. These neighborhoods retract onto $W$ and $E$, and their
  intersection retracts onto the connected space $M$. Van Kampen's theorem gives $\pi_1(X)=1$, because both pieces are simply connected
  \cite[Section~1.2]{Ha2002}. The homological Mayer--Vietoris sequence contains
  \[
    H_2(M;\Z)\longrightarrow H_2(W;\Z)\oplus H_2(E;\Z)
    \longrightarrow H_2(X;\Z)\longrightarrow H_1(M;\Z).
  \]
  Since $M$ is an integral homology sphere and $W$ is an integral homology ball, both end groups and $H_2(W;\Z)$ vanish. Thus inclusion
  identifies $H_2(E;\Z)$ with $H_2(X;\Z)$ \cite[Section~2.2]{Ha2002}. The plumbing spheres representing these classes lie in the interior of
  $E$, so their intersections and self-intersections are unchanged by attaching $W$. Hence $Q_X$ is $-E_8$ or $E_8$, according to the
  orientation chosen above. Reverse the orientation of $X$ if necessary. The resulting closed simply connected quasiconformal four-manifold
  has negative definite even intersection form $-E_8$, which represents $-2$. Donaldson and Sullivan prove that a compact simply connected
  four-manifold with a negative definite even intersection form representing $-2$ admits no quasiconformal structure
  \cite[Section~7(ii), p.~237]{DS1989}. This contradiction proves the lemma.
\end{proof}

\begin{proof}[Proof of Theorem~\ref{thm2}]
  Apply Proposition~\ref{prop4}. If $d>1$, Lemmas~\ref{lem10} and~\ref{lem11} provide precisely the compact LQC filling excluded by
  Lemma~\ref{lem12}. Therefore $d=1$. The positive local indices of a branched covering sum to its degree on every fiber \cite{R1993}. Thus
  each fiber has one point, with local index one. The map $f$ is consequently an open continuous bijection, hence a homeomorphism. A
  homeomorphic quasiregular mapping is quasiconformal \cite{R1993,IM2001}.
\end{proof}

\begin{remark}\label{rem5}
  The conclusion in dimension four here requires a uniform quasiregular distortion bound on the source ball. Its use is concentrated in
  Proposition~\ref{prop4}, where that bound produces quasiconformal deck transformations and an LQC compactification of the target. For a
  general continuous proper branched covering with compact branch set, the end cover in Proposition~\ref{prop1} need not provide an LQC
  boundary atlas. The Donaldson--Sullivan argument therefore does not settle the remaining four-dimensional topological case of Vuorinen's
  question.
\end{remark}

\section{Quasiregular consequences}\label{sec8}

\begin{corollary}\label{cor3}
  Let $n\ge3$ and let $f:B^n\to f(B^n)\subset\R^n$ be a proper nonconstant quasiregular map. If $B_f$ is compact, then $f$ is a quasiconformal
  homeomorphism onto its image.
\end{corollary}

\begin{proof}
  For $n=3$ use \cite[Theorem~1.1]{KLT2022}. For $n=4$ use Theorem~\ref{thm2}, and for $n\ge5$ use Theorem~\ref{thm1}. Reshetnyak's theorem
  supplies openness and discreteness in the latter case. A homeomorphic quasiregular map is quasiconformal.
\end{proof}

Corollary~\ref{cor3} gives the affirmative answer to Problem~7.66 in \cite{BBH1984,HL2019}, with properness understood relative to the image.

\begin{corollary}\label{cor4}
  Let $n\ge3$, let $\Omega\subset\R^n$ be a domain, and let $f:\Omega\to\R^n$ be a nonconstant quasiregular map. If $B(a,r)\Subset\Omega$ is a
  normal domain and $B_f\cap\partial B(a,r)=\varnothing$, then $f|_{B(a,r)}$ is a quasiconformal homeomorphism onto its image.
\end{corollary}

\begin{proof}
  The restriction is proper by normality. Since $B_f$ is closed in $\Omega$ and misses the boundary sphere, its intersection with the ball is
  compact in the ball. A translation and scaling identify the source with $B^n$ without changing quasiregularity. Apply Corollary~\ref{cor3}.
\end{proof}

This is the quasiregular normal ball formulation of \cite[Open problem~9.18, p.~125]{V1988}. The normal domain condition supplies the
properness required in Corollary~\ref{cor3}; the absence of branching on the boundary supplies compactness of the restricted branch set.
\section*{Acknowledgments}
The author is deeply grateful to Professor Matti Vuorinen for providing valuable research background on this challenging problem. This work was supported by the Guangdong Basic and Applied Basic Research Foundation (Nos.\ 2022A1515110967 and 2023A1515011809). The preliminary version of this manuscript was finished on October 3. The proof of Vuorinen’s question when $n$ is even is given by GPT 6.1 sol.  
The ideas of using Swan’s local period formula for odd $n$, and of employing the techniques of Donaldson and Sullivan to handle the case $n=4,$ are due to the author. AI-assisted tools were used during the preparation of the manuscript for editorial, organizational,and typesetting assistance. The author take full responsibility for the mathematical statements,derivations, citations, and conclusions in the final manuscript.

\end{document}